\documentclass[12pt]{article}

\usepackage[dvips]{graphicx}
\usepackage{t1enc,amssymb,amsbsy,amsthm,amsmath,amsfonts}

\def\th{\hbox {\rm tanh\,}}

\def\ch{\hbox {\rm cosh\,}}
\def\cth{\hbox {\rm ctnh\,}}
\def\sh{\hbox {\rm sinh\,}}
\def\tg{\hbox {\rm tan\,}}
\def\ctg{\hbox {\rm ctn\,}}
\def\arctg{\hbox {\rm arctan\,}}
\def\arcsin{\hbox {\rm arcsin\,}}
\def\ln{\hbox {\rm ln\,}}
\def\Arth{\hbox {\rm artanh\,}}
\def\Arsh{\hbox {\rm arsinh\,}}

\renewcommand{\Re}{{\rm I\kern-0.16em R}}

\def\@begintheorem#1#2{\trivlist \item[\hskip \labelsep{\bf #1\ #2}]}
\def\@opargbegintheorem#1#2#3{\trivlist
      \item[\hskip \labelsep{\bf #1\ #2\ (#3)}]}

\newtheorem{proposition}{Proposition}[section] %denna med i latex

\newtheorem{theorem}[proposition]{Theorem}

\newtheorem{example}[proposition]{Example}
\newtheorem{remark}[proposition]{Remark}

\def\ch{\hbox {\rm cosh}}
\def\sh{\hbox {\rm sinh}}
\def\P{{\bf P}}
\def\R{{\bf R}}

\def\E{{\bf E}}

\def\cF{{\cal F}}

\def\cG{{\cal G}}

\def\al{\alpha}

\def\Ga{\Gamma}
\def\ga{\gamma}
\def\si{\sigma}
\def\te{\theta}
\def\be{\beta}
\def\la{\lambda}
\def\om{\omega}

\numberwithin{equation}{section}

\begin{document}

\author{
Andrei Borodin\\
{\small Steklov Mathematical Institute,}\\
{\small St. Petersburg Division,}\\
{\small Fontanka 27}\\
{\small 191011 St. Petersburg, Russia}\\
{\small email: borodin@pdmi.ras.ru}
\and
Paavo Salminen\\{\small Åbo Akademi,}\\
{\small Mathematical Department,}\\
{\small V\"anriksgatan 3 B}\\
{\small FIN-20500 Åbo, Finland,} \\
{\small email: phsalmin@abo.fi}
}
\vskip5cm

\title{On some exponential integral functionals of BM($\mu$) and
  BES(3)}
%rownian motion with drift}

\maketitle

\begin{abstract}
In this paper we derive the Laplace transforms of the integral
functionals 
$$
\int_0^\infty 
\left(p\left(\exp(B^{(\mu)}_t)+1\right)^{-1}+
q\left(\exp(B^{(\mu)}_t)+1\right)^{-2}\right)\, dt,
$$
%$$
%\int_0^\infty 
%\Big(\frac{p}{\exp(B^{(\mu)}_t)+1}+
%\frac{q}{\big(\exp(B^{(\mu)}_t)+1\big)^{2}}\Big)\, dt,
%$$
and
%$$
%\int_0^\infty 
%\Big(\frac{p}{\exp(R^{(3)}_t)-1}+
%\frac{q}{\big(\exp(R^{(3)}_t)-1\big)^{2}}\Big)\, dt,
%$$
$$
\int_0^\infty 
\left(p\left(\exp(R^{(3)}_t)-1\right)^{-1}+
q\left(\exp(R^{(3)}_t)-1\right)^{-2}\right)\, dt,
$$
where $p$ and $q$ are real numbers, $\{B^{(\mu)}_t:\ t\geq 0\}$ 
is a Brownian motion with drift 
$\mu>0,$ BM($\mu$), and $\{R^{(3)}_t:\ t\geq 0\}$ is a 
3-dimensional Bessel process, BES(3).
The transforms are given in terms of Gauss' hypergeometric
functions and it is seen that the results are closely related to some
functionals of Jacobi diffusions.
This work generalizes and completes some results of Donati--Martin and
Yor \cite{donatimartinyor97} and Salminen and Yor \cite{salminenyor04}.
%\cite{salminenyor04} the functional with $p=0$ and $q=1$ is 
%studied and explicit analysis is done for the case $\mu=1/2.$
%The present work extends some of these results.
\\ \\%\bigskip\noindent
\\ \\ %\bigskip\noindent
\end{abstract}

\section{Introduction}
\label{sec0}
Let $B^{(\mu)}=\{ B^{(\mu)}_t:=B_t+\mu t \,:\, t\geq
0\}$ be a Brownian motion with drift $\mu>0$ started from 0. For 
the
non-negative and locally integrable function $f$ 
consider the functional 
$$
I_\infty(f):=\int_0^\infty f(B^{(\mu)}_s)\,ds.
$$
It is known, see Engelbert and Senf \cite{engelbertsenf91} and Salminen
and Yor \cite{salminenyor03}, 
\begin{equation}
\label{e01}
I_\infty(f)<\infty\ {\rm a.s}\quad \Leftrightarrow\quad \int^{\infty}
f(x)\, dx<\infty.
\end{equation}
In \cite{salminenyor03} also other properties of $I_\infty(f)$ are
studied; in particular, existence of moments.

During the recent years and because of 
applications in financial mathematics
much interest has been focused on 
exponential functionals of $B^{(\mu)}.$ 
We refer to Yor \cite{yor01} for a collection of papers on this topic 
and also for further references.  
%see, e.g.,
% Dufresne \cite{dufresne90}, 
%\cite{dufresne01}, 
%Yor \cite{yor92a}, \cite{yor92b}, \cite{yor01},  
%Donati--Martin and Yor \cite{donatimartinyor97},
%Donati-Martin, Matsumoto and Yor \cite{donatimatsumotoyor00a},
%\cite{donatimatsumotoyor00b}, \cite{donatimatsumotoyor00c},
%\cite{donatimatsumotoyor00d},
%Matsumoto and Yor \cite{matsumotoyor98}, \cite{matsumotoyor99a}, 
%\cite{matsumotoyor99b}, \cite{matsumotoyor00}, 
%\cite{matsumotoyor01}, \cite{matsumotoyor03},
%Salminen and Yor \cite{salminenyor03}, \cite{salminenyor03b},
%\cite{salminenyor04}. 
A particular rôle in these investigations has been played by 
the functional   
\begin{equation}
\label{du} 
\int_0^\infty  \exp(-2B_s^{(\mu)})\, ds
\end{equation}
which, as shown in Dufresne \cite{dufresne90}, can be obtained 
as a perpetuity in a discrete pension funding 
scheme after a limiting procedure. See also Yor \cite{yor92b}.

In \cite{salminenyor03b} and \cite{salminenyor04} the functional   
\begin{equation}
\label{y}
\int_0^{\infty}  (a+\exp(B^{(\mu)}_s))^{-2}\, ds
\end{equation}
is analyzed the motivation being that the functional in (\ref{y}) 
has all the moments (in fact, some exponential moments) when 
Dufresne's functional (\ref{du}) has only some moments. 
%For a more detailed discussion and further 
%references we refer to \cite{salminenyor03b} and 
The distribution of the functional in 
(\ref{y}) is computed in \cite{salminenyor03b} for $\mu=1/2$ on one hand by 
connecting the functional to a hitting time of another Brownian motion
with drift via random time change and on the other hand 
using the Feynman-Kac formula. 

In this paper we derive %using the Feynman-Kac formula 
the Laplace transform of
the functional in (\ref{y}) for arbitrary $\mu>0.$ 
In fact, we study the more general functional 
\begin{equation}
\label{a}
\int_0^\infty
\Big(\,
p\,\big(\exp(B^{(\mu)}_t)+1\big)^{-1} +
q\,\big(\exp(B^{(\mu)}_t)+1\big)^{-2}\Big)\, dt.
\end{equation}
Notice that choosing here, e.g. , $p=-q=1$
we obtain the functional
\begin{equation}
\label{v}
\int_0^\infty
\cosh^{-2}(B^{(\mu)}_t/2)\, dt.
\end{equation}
The key observation is that the functional in (\ref{a}) can be 
expressed, roughly speaking, via the first hitting time of 
a Jacobi diffusion the Laplace transform of which can be 
computed explicitly in terms of Gauss' hypergeometric functions. It is 
also seen that similar techniques apply to compute the Laplace transform of 
an analogous functional for three-dimensional Bessel process 
$\{R_t\, ; t\geq 0\}$:
\begin{equation}
\label{aa}
\int_0^\infty 
\left(p\left(\exp(R_t-1\right)^{-1}+
q\left(\exp(R_t)-1\right)^{-2}\right)\, dt.
\end{equation}

The paper is organized so that in the next section some general results on 
the random time techniques connecting integral functionals and hitting
times are presented. In this section we also discuss the Feynman--Kac
approach. Definitions and some properties of Gauss' hypergeometric
functions are given in Section \ref{sec1}. In Section \ref{sec2} we
discuss some properties of Jacobi diffusions focusing on two special cases  
emerging from our computations. The functional in (\ref{a}) is treated in
Section \ref{sec3} and (\ref{aa}) in
Section \ref{sec4}. 

%h

\section{Transforming diffusions by changing time and space}
\label{sec11}

In this section we give a fairly straightforward generalization, needed in Section
\ref{sec2}, of the main result in Salminen and Yor
\cite{salminenyor04} for quite general
diffusions. We also discuss the corresponding Feynman--Kac approach.   

Let $B$ be a Brownian motion and
consider a diffusion $Y$ on an open interval $I=(l,r)$
determined by the SDE
\begin{equation}
\label{a1}
dY_t=\sigma(Y_t)\,dB_t+b(Y_t)\,dt.% \qquad Y_0=x\in I.
\end{equation}
It is assumed that $\sigma$ and $b$ are continuous and
$\si(x)>0$ for $x\in I$. The diffusion $Y$ is considered up to
its explosion (or life) time
$$
\zeta:=\inf\{t\ :\ Y_t\not\in I\}.
$$

\begin{proposition}
\label{!}
{\sl
Let $g(x)$, $x\in I$, be a twice 
continuously differentiable function such that $g'(x)\not=0$. 
Consider the integral functional
$$
A_t:=\int_0^t \big(g'(Y_s)\si(Y_s)\big)^2\,ds, \qquad t\in [0,\zeta),
$$
and its inverse
$$
a_t:=\min\big\{s: A_s>t\big\}, \qquad t\in [0,A_\zeta).
$$
Then the process $Z$ given by
\begin{equation}
\label{a2}
Z_t:=g\left(Y_{a_t}\right), \qquad t\in[0,A_\zeta),
\end{equation}
is a diffusion satisfying the SDE
\begin{equation}
\label{a3} 
dZ_t=d\widetilde B_t+G(g^{-1}(Z_t))\,dt, \qquad t\in[0,A_\zeta).
\end{equation}
where $\widetilde B_t$ is a Brownian motion and
\begin{equation}
\label{a4} 
G(x)=\left(g'(x)\si(x)\right)^{-2}\left( \frac 12\, \si(x)^2\, g''(x)
+b(x)\,g'(x)\right).
\end{equation}
}
\end{proposition}

\begin{proof} By Ito's formula for $u\leq \zeta$
$$
g(Y_u)-g(Y_0)=\int_0^u g'(Y_s)\si(Y_s)\,dB_s
+\int_0^u\big(2^{-1}g''(Y_s)\si^2(Y_s)+g'(Y_s)b(Y_s)\big)\,ds.
$$
Replacing $u$ by $a_t$ yields
$$
Z_t-Z_0=\int_0^{a_t}g'(Y_s)\si(Y_s)\,dB_s
+\int_0^{a_t}\big(g'(Y_s)\si(Y_s)\big)^2G(Y_s)\,ds.
$$
Since $a_t$ is the inverse of $A_t$ and $A_s'=\big(g'(Y_s)\sigma(Y_s)\big)^2$
we have
\begin{equation}
\label{a5}
a_t'=\frac{1}{A'_{a_t}}
=\left(g'(Y_{a_t})\sigma(Y_{a_t})\right)^{-2}.
\end{equation}
From L\'evy's theorem it follows that
$$
\widetilde B_t:=\int_0^{a_t}g'(Y_s)\sigma(Y_s)\,dB_s,
\qquad t\in[0,A_\zeta).
$$
is a (stopped) Brownian motion. As a result we obtain for $t<A_\eta$
\begin{eqnarray*}
&&
Z_t-Z_0=\widetilde B_t+\int_0^t
\big(g'(Y_{a_s})\si(Y_{a_s})\big)^2G(Y_{a_s})\,da_s\\
&&\hskip1.5cm
=\widetilde B_t+\int_0^t G(g^{-1}(Z_s))\,ds.
\end{eqnarray*}
\end{proof}

It is interesting and useful to understand how the relation 
(\ref{a2}) expresses itself on the level of differential 
equations associated with the 
diffusions $Y$ and $Z.$ In particular, we consider in this framework 
the Feynman--Kac formula. First, we give the following proposition which is a 
generalization of the result for Brownian motion stated in 
\cite{borodinsalminen02} VI.15 p.113 and concerns  distributions of
functionals of Brownian motion stopped according to the
inverse of an additive functional.

\begin{proposition}
\label{fk}
{\sl
Let $F$ and $f$ be continuous functions on $I$. Assume also that  $F$
is bounded and $f$ is non-negative. Let $\tau$ be an exponentially
(with parameter $\lambda$) distributed random variable independent of
$Y.$ Then the function
$$
U(x):=\E_x\left(F(Y_{a_\tau})\exp\left(-\int_0^{a_\tau}f(Y_s)\,ds
\right)\right)
$$
is the unique bounded solution of the differential equation
\begin{eqnarray}
\label{a7}
\nonumber
&&\hskip-1cm
\frac{\si^2(x)}2 u''(x)+b(x)u'(x)-\big(\la\big(g'(x)\si(x)\big)^2+f(x)\big)u(x)
\\
&&\hskip5cm
=-\la\big(g'(x)\si(x)\big)^2 F(x).
\end{eqnarray}
}
\end{proposition}
\noindent
Next, introduce for $x\in J:=g(I)$
\begin{equation}
\label{a8}
Q(x):=U(g^{-1}(x)).
\end{equation}
Straightforward calculations show that the function $Q$
satisfies
\begin{eqnarray}
\label{a9}
\nonumber
&&\hskip-1cm
\frac12 Q''(x)+G(g^{-1}(x))Q'(x)-\Big(\la+
\frac{f(g^{-1}(x))}{\big(g'(g^{-1}(x))\si(g^{-1}(x))\big)^2}\Big)Q(x)
\\
&&\hskip5cm
=-\la F(g^{-1}(x)).
\end{eqnarray}
By the classical Feynman--Kac result the equation (\ref{a9}) has
only one bounded and non-negative solution, and, consequently, we have the
following probabilistic interpretation
$$
Q(x)=\E_x\left (F(g^{-1}(Z^\circ_\tau))\exp\Big(-\int_0^\tau
\frac{f(g^{-1}(Z^\circ_s))}
{\big(g'(g^{-1}(Z^\circ_s))\si(g^{-1}(Z^\circ_s))\big)^2}\,ds\Big)\right ),
$$
where the diffusion $Z^\circ$ has the generator
$$
\frac12 \frac{d^2}{dx^2}+G(g^{-1}(x))\frac{d}{dx},
$$
and can be viewed as a solution of the SDE (\ref{a3}).
Using (\ref{a8}) and the probabilistic interpretations of $U$ and $Q$
we obtain
\begin{eqnarray}
\label{i2}
%&&
%\E_{g^{-1}(x)}\left (F(Y_{a_\tau})
%\exp\Big(-\int_0^{a_\tau} f(Y_{a_s})\,ds\Big)\right )
%\\
%&&
%\hskip1cm
%=\E_{g^{-1}(x)}\left (F(Y_{a_\tau})
%\exp\Big(-\int_0^{\tau} f(Y_{a_s})\,da_s\Big)\right )
%\\
%
&&
\E_{g^{-1}(x)}\left(F(Y_{a_\tau})
\exp\Big(-\int_0^\tau \frac{f(Y_{a_s})}
{\big(g'(Y_{a_s})\si(Y_{a_s})\big)^2}\,ds\Big)\right)
\\
&&
\nonumber
\hskip1cm
=
\E_x\left (F(g^{-1}(Z^\circ_\tau))\exp\Big(-\int_0^\tau
\frac{f(g^{-1}(Z^\circ_s))}
{\big(g'(g^{-1}(Z^\circ_s))\si(g^{-1}(Z^\circ_s))\big)^2}\,ds\Big)\right ).
\end{eqnarray}
From this identity we can recover, e.g., that the distributions
of $Y_{a_t}$ and $Z^\circ_t$ are equal (from (\ref{a2}) we, of course, 
know much more). In this paper, see Sections
4.2 and 4.3, the
identity (\ref{i2}) is used in a particular case to find the solutions
of the equation (\ref{a9}) in terms of the (known) solutions of the
equation (\ref{a7}). Notice also that Proposition \ref{fk} can be proved
by using first the classical Feynman--Kac formula for (\ref{a9}) and
then Proposition \ref{!}.

Next two propositions connect occupation time functionals 
with first hitting times.
We use the notation in Proposition \ref{!}.

\begin{proposition}
\label{!!}
{\sl
{\bf 1.} Let $x\in I$ and $y\in I$
be such that $P_x$-a.s.
$$
H_y(Y):=\inf\{t:\, Y_t=y\}<\infty.
$$
Then, under the assumptions in Proposition \ref{!},
\begin{equation}
\label{aa10}
A_{H_y(Y)}=\inf\{t:\, Z_t=g(y)\}=: H_{g(y)}(Z)\quad a.s.
\end{equation}
with $Y_0=x$ and $Z_0=g(x).$

\noindent
{\bf 2.} Assume $g(r):=\lim_{z\to r}g(z)$ exists 
(recall that $r$ denotes the right hand side end 
point of $I$). Suppose also that each of the following statements holds a.s.
\begin{equation}
\label{aa100}
(i)\ \
\zeta=+\infty,
\quad 
(ii)\ \ \lim_{t\to\infty}Y_t=r,\quad 
(iii)\ \ A_\infty:=\lim_{t\to +\infty}A_t<\infty .
\end{equation}
%\begin{equation}
%\label{aa100}
%A_\infty:=\lim_{t\to +\infty}A_t<\infty .
%\end{equation}
Then 
\begin{equation}
\label{a10}
A_\infty= H_{g(r)}(Z)\quad a.s.%\inf\{t:\, Z_t=g(r)\}=:
\end{equation}
}
\end{proposition}
 
\begin{proof}
To prove (\ref{aa10}), notice that from the definition of $Z$ (see (\ref{a2}))
we obtain 
$$
Z_{A_t}=g(Y_t).
$$
Letting here $t\to H_y(Y)$ and using the fact that $g$ is monotone yield 
(\ref{aa10}). The claim (\ref{a10}) is proved similarly. 

\end{proof}

\begin{remark}
{\rm 
 A sufficient condition for (iii) in (\ref{aa100}) is clearly that 
the mean of $A_\infty$ is finite:
\begin{eqnarray*}
&&\E_x\left(A_\infty\right)= 
\int_0^\infty \E_x\left(\left(g'(Y_s)\si(Y_s)\right)^2\right)\,ds \\
&&\hskip1.7cm
= \int_l^r G_0(x,y)\, (g'(y)\si(y))^2\,m(dx)<\infty,
\end{eqnarray*}
where $G_0$ denotes the Green kernel of $Y$ and $m$ is the speed measure
(for these see, e.g., Borodin and Salminen \cite{borodinsalminen02}).
For a neccessary and sufficient condition in the case of a 
Brownian motion with drift, see (\ref{e01}).
%Engelbert and Senf \cite{engelbertsenf} and 
%Salminen and Yor \cite{salminenyor03}).
}
\end{remark}

For the next proposition recall from (\ref{a5}) and (\ref{a2}) 
that the inverse $a_t$ 
of $A_t$ is given for $t<A_\zeta$ by 
\begin{equation}
\label{aa6}
a_t=\int_0^t \left(g'(g^{-1}(Z_s))\,\sigma(g^{-1}(Z_s))\right)^{-2}\,ds.
\end{equation}

\begin{proposition}
\label{!!!}
Assume that each of the following statements holds a.s.  
$$(i)\ \  \zeta<\infty,\quad (ii)\ \  \lim_{t\to\zeta}Y_t=l,\quad 
(iii)\ \ A_{\zeta}=\infty.
$$ 
Then 
\begin{equation}
\label{a6}
a_\infty:=\int_0^\infty \left(g'(g^{-1}(Z_s))\,
\sigma(g^{-1}(Z_s))\right)^{-2}\,ds
=\zeta\quad {\rm a.s.}
%H^Y_l\quad {\rm a.s.}
\end{equation}
\end{proposition}

\begin{proof} The claim is immediate from (\ref{aa6}) using 
the assumptions (i)--(iii) and the fact that 
$a_t$ is the inverse of $A_t.$ Notice also that in this case 
$\zeta=\lim_{x\to l}H_x(Y).$ 
\end{proof}

We conclude this section by a proposition which relates the results in 
Propositions \ref{!!} and \ref{!!!} to the Feynman--Kac formula. We consider 
only the case in Proposition \ref{!!}.2 and leave the other cases to the 
reader. 

\begin{proposition}
\label{rk!}
{\bf 1.} Under the assumptions of Proposition \ref{!!}.2,  the function
$$
\Phi(x):=\E_x\left(\exp\left(-\rho\int_0^\infty 
\big(g'(Y_s)\si(Y_s)\big)^2\,ds\right)\right), \quad x\in I,\quad \rho>0
$$
is the unique increasing positive solution of the differential equation  
\begin{equation}
\label{a11} 
\frac 12\,\si^2(x)\, \Phi''(x)+b(x)\,\Phi'(x)
-\rho\,\big(g'(x)\si(x)\big)^2\,\Phi(x)=0,
\end{equation}
satisfying $\lim_{x\to r}\Phi(x)=1$.

\noindent
{\bf 2.} Define for $x\in g^{-1}(I)$
%\begin{equation}
%\label{a12}
$$
\Psi(x):=\Phi(g^{-1}(x)).
$$
%\end{equation}
Then $\Psi$ is the unique increasing positive solution 
of the differential equation
$$
\frac 12\, \Psi''(x)+G(g^{-1}(x))\,\Psi'(x)-\rho\, \Psi(x)=0
$$
satisfying $\lim_{x\to r}\Psi(x)=1$. Moreover, 
\begin{eqnarray}
\label{a12}
&&
\nonumber
\Psi(x)=\E_x\left(\exp(-\rho H_{g(r)}(Z))\right)\\
&&
\hskip1cm
=
\E_{g^{-1}(x)}\left(\exp\left(-\rho\int_0^\infty 
\left(g'(Y_s)\si(Y_s)\right)^2\,ds\right)
\right).
\end{eqnarray}
\end{proposition}

\begin{proof}
The first part of the proposition can be proved as Proposition XX in Salminen 
and Yor \cite{salminenyor04}. The fact that $\Psi$ has the claimed properties 
is an easy exercise in differentiation. The first equality in (\ref{a12}) is 
standard diffusion theory and the second one follows from the 
definition of $\Psi.$
\end{proof}

\section{Gauss' hypergeometric functions}
\label{sec1}
We start with by recalling from Abramowitz and Stegun
\cite{abramowitzstegun70} p. 556,
the definition of Gauss' hypergeometric series (or functions):
\begin{equation}
\label{F}
F(\alpha,\beta,\gamma;x):=
\frac{\Gamma(\gamma)}{\Gamma(\beta)\Gamma(\alpha)}
\sum_{n=0}^\infty\frac{\Gamma(\alpha+n)\Gamma(\beta+n)}
{\Gamma(\gamma+n)}\, \frac{x^n}{n!}.
\end{equation}
It is easily seen that $F$ is well defined for $|x|<1$
but we consider $F$ only for real values on $x$ such that $0<x<1.$
However, we allow complex conjugate values for $\alpha,$ and $\beta,$
but take $\gamma$ real (see Theorems \ref{main} and \ref{main2}).
Notice that
$$
F(\alpha,\beta,\gamma;x)=F(\beta,\alpha,\gamma;x).
$$
%where $\alpha>0$ and $\alpha +1>\gamma>\beta>0.$
It is well known (see Abramowitz and Stegun \cite{abramowitzstegun70}
or Lebedev \cite{lebedev72} p. 162)  that $F$ is a solution of
the ODE
\begin{equation}
\label{gauss}
x(1-x)\, v''(x)+(\gamma -(\alpha+\beta+1)\, x)\,v'(x)-\alpha\,\beta \,v(x)=0.
\end{equation}
Straightforward calculations show that also
\begin{equation}
\label{G}
\widehat F(\alpha,\beta,\gamma;x):=F(\alpha,\beta,\alpha+\beta+1-\gamma;1-x)
\end{equation}
is a solution of (\ref{gauss}).
%In particular, if
%\begin{equation}
%\label{c1}
%\alpha,\beta,\gamma\,\geq 0,
%\end{equation}
%it
%follows from (\ref{F}) that  $F$ is increasing function of $x$
%and, therefore, $\widehat F$ is decreasing if
%\begin{equation}
%\label{c2}
%\gamma<\alpha+\beta+1.
%\end{equation}
%Consequently, we have found under
%the conditions (\ref{c1}) and (\ref{c2})
%two linearly independent solutions of
%(\ref{gauss}) which generate all the solutions of
%(\ref{gauss}).

Recall also that for 
$0<x<1,$ and ${\rm Re}(\ga)>{\rm Re}(\alpha)>0$
it holds 
\begin{equation}
\label{F2}
F(\alpha,\beta,\gamma;x)=
\frac{\Gamma(\gamma)}{\Gamma(\al)\Gamma(\gamma-\al)}
\int_0^1 t^{\al-1}(1-t)^{\gamma-\al-1}(1-tx)^{-\beta}\, dt.
\end{equation}
From (\ref{F2}) we obtain by changing variables 
\begin{eqnarray}
\label{G2}
&&
\nonumber
\hskip-1.3cm
\widehat F(\alpha,\beta,\gamma;x)=F(\alpha,\beta,\alpha+\beta+1-\gamma;1-x)\\
&&\hskip1cm
=\frac{\Gamma(\alpha+\beta+1-\gamma)}{\Gamma(\al)\Gamma(\beta+1-\gamma)}
%&&\hskip4cm
%\nonumber
%\times
\int_0^\infty t^{\al-1}(1+t)^{\gamma-\al-1}(1+tx)^{-\beta}\, dt,
\end{eqnarray}
which is valid for $0<x<1,$ 
${\rm Re}(\be+1-\ga)>0,$ and ${\rm Re}(\al)>0.$ 
Notice that 
$$
\widehat F(\alpha,\beta,\gamma;1)=1
$$
and 
\begin{equation}
\label{Glimit}
\widehat F(\alpha,\beta,\gamma;0)=
\frac{\Gamma(\alpha+\beta+1-\gamma)\,\Gamma(1-\ga)}
{\Gamma(\al+1-\ga)\,\Gamma(\beta+1-\gamma)},
\end{equation}
where all the arguments in our case (see Theorem \ref{main})
are positive or have positive real parts.

%From the the integral representation (\ref{F2}) it is seen by monotone
%convergence for $\gamma>\alpha+\beta$
%\begin{equation}
%\label{limit}
%F(\alpha,\beta,\gamma;1):=\lim_{x\to 1-}F(\alpha,\beta,\gamma;x)
%=\frac{\Gamma(\gamma)\,\Gamma(\gamma-\alpha-\beta)}
%{\Gamma(\gamma-\alpha)\Gamma(\gamma-\beta)}.
%\end{equation}

\section{Jacobi diffusion}
\label{sec2}
\subsection{Definition and two particular cases}
\label{sec21}

Consider the diffusion  $X=\{ X_t\,:\, t\geq 0\}$ living in $(0,1)$
satisfying the SDE
\begin{equation}
\label{SDE}
dX_t=\sqrt{2X_t(1-X_t)}\,dB_t+(\gamma-(\alpha+\beta +1)X_t)\,dt.
\end{equation}
Following
Karlin and Taylor \cite{karlintaylor81} p. 335, $X$ is called 
a Jacobi diffusion with parameters  $\al,$ $\be,$ and $\ga.$ See 
the papers by  Mazet
\cite{mazet97},
 Warren and Yor \cite{warrenyor98}, Hu, Shi and Yor \cite{hushiyor99},
Schoutens \cite{schoutens00}, and
Albanese and Kuznetsov \cite{albanesekuznetsov03} for some results and
applications of Jacobi diffusions (and also for further references).
Next proposition gives a basic simple property of Jacobi diffusions.

\begin{proposition}
\label{1-x}
{\sl
Let $X$ be a Jacobi diffusion with parameters $\al,$ $\be,$ and $\ga.$ 
Then $Y:=1-X$ is a Jacobi diffusion with parameters $\al,$ 
$\be,$ and $\al+\be+1-\ga.$
}
\end{proposition}

The generator of $X$ is given by
$$
\cG u(x)=
x(1-x)\,u''(x)+(\gamma-(\alpha+\beta +1)\,x)u'(x),
$$
and we use
$$
S(x)=\int_{1/2}^x y^{-\gamma}(1-y)^{\gamma-\alpha-\beta -1}\, dy
$$
as the scale, and
$$
m(x_1,x_2)=\int_{x_1}^{x_2} y^{\gamma-1}(1-y)^{\al+\be-\ga}
\, dy,\quad  0\leq x_1<x_2\leq 1
$$
as the speed.
Obviously, $X$ exhibits very different behaviour
when the values of the
parameters $\al,$ $\be,$ and $\ga.$ are varied.
In view of our applications we consider two cases.
\hfill\break
\vskip.1cm\noindent
{\bf Case 1}: $\gamma<1$ and $\alpha+\beta>1.$
\hfill\break
\vskip.1cm
Here, $S(0+)>-\infty$ and $S(1-)=+\infty,$ and
by the standard boundary point analysis (see, e.g., Borodin and
Salminen \cite{borodinsalminen02})
\begin{itemize}
\item{} if $\ga\leq 0$ then 0 is exit-entrance (regular),
\item{} if $0<\ga<1$ then 0 is exit-not-entrance,
\item{} 1 is entrance-not-exit.
\end{itemize}
Consequently, $X$ is transient and
$$
H_0(X):=\inf\{t\,:\, X_t=0\}<\infty\quad  {\rm a.s.}
$$
In the regular case, i.e., $\ga\leq 0$, we choose $0$ to be a killing
boundary.
We remark that the boundaries 0 and 1 are as displayed above also if the
condition $\alpha+\beta>1$ is extended to
$\alpha+\beta>\ga-1.$
\hfill\break
\vskip.1cm\noindent
{\bf Case 2}: $\gamma=\alpha+\beta\geq 2.$
\hfill\break
\vskip.1mm
Now we have $S(0+)=-\infty$ and $S(1-)=+\infty,$
$ m(0,1)<\infty,$  and
\begin{itemize}
\item{} 0 is entrance-not-exit,
\item{} 1 is entrance-not-exit.
\end{itemize}
Hence, $X$ is positively recurrent. Notice, cf. Proposition \ref{1-x},
that if $X^{(1)}$
satisfies (\ref{SDE}) with $\ga=\alpha+\beta$ then
$X^{(2)}:=1-X^{(1)}$
solves (\ref{SDE}) with $\ga=1.$ This fact
is also transparent in the expressions for $S$ and $m.$

\subsection{Functionals of Jacobi diffusion; Case 1}
\label{sec22}

Consider the Jacobi diffusion given by (\ref{SDE}) with $\ga<1$ and
$\al+\be >1.$ Let for $c,\,\te\geq 0$
$$
f(x):=\frac{1}{\te}\log\left(\frac{1-x}{c\,x}\right),\quad 0<x<1.
$$
Clearly, $f$ is decreasing, $f(0+)=+\infty,$ 
$f(1-)=-\infty$, and
$$
f^{-1}(x)=\frac{1}{c\,{\rm e}^{\,\te x}+1}.
$$
We are interested in the process $\{f(X_t)\,:\,t\geq 0\}.$ 
By Ito's formula for $t<H_0(X)$
\begin{eqnarray}
\label{e1}
\nonumber
&&
\log\left(\frac{1-X_t}{c\,X_t}\right)-
\log\left(\frac{1-X_0}{c \,X_0}\right)\\
&&\hskip.5cm
=\int_0^t \frac{-\sqrt{2}}{\sqrt{X_s(1-X_s)}}\,dB_s
+ \int_0^t \frac{1-\ga+(\al+\be-1)X_s}{X_s(1-X_s)}\,ds.
\end{eqnarray}
Define for $t<H_0(X)$
$$
A_t:=\frac{2}{\te^2}\int_0^t \frac{1}{X_s(1-X_s)}\,ds,
$$
and set $A_t=+\infty$ for $t\geq H_0(X)$. Let $a$ be the inverse of
$A,$ i.e., 
$$
a_t:=\inf\{s\,:\, A_s>t\},
$$ 
and notice that $a_t\leq
H_0(X)$ for all $t\geq 0.$
Moreover,
\begin{equation}
\label{alpha1}
a_t=\frac{\te^2}{2}\int_0^t X_{a_s}(1-X_{a_s})\,ds.
\end{equation}
The process given by
\begin{equation}
\label{Z}
Z_t:=\frac{1}{\te}\log\left(\frac{1-X_{a_t}}{c\,X_{a_t}}\right)
\end{equation}
is well defined for all $t$ such that $a_t<H_0(X),$
and from (\ref{e1}) we obtain
\begin{equation}
\label{ZZ}
Z_t-Z_0= B^\circ_t+\frac {\te}{ 2}\int_0^t\left(1-\ga+\frac{\al+\be-1}
{c{\rm e}^{\,\te Z_s}+1}\right)\,ds,
\end{equation}
where $B^\circ $ is a Brownian motion. From (\ref{ZZ}) it is seen that
$Z$ is, in fact, non-exploading and, therefore, from (\ref{Z}) it
follows that $a_t<H_0(X)$ for all $t$ and $a_t\to H_0(X)$ as
$t\to\infty.$ Hence, it also holds that $A_t\to\infty$ as $t\to
H_0(X),$ and we are in the case treated in Proposition \ref{!!!}.
We have 
\begin{eqnarray*}
&&
\int_0^t\frac{c\,\exp(\,\te Z_s)}{(c\,\exp(\te Z_s)+1)^2}\, ds =
\int_0^t \frac{1-X_{a_s}}{X_{a_s}}
\left(\frac{1-X_{a_s}}{X_{a_s}}+1\right)^{-2}\, ds\\
&&\hskip4cm
=\int_0^t X_{a_s}(1-X_{a_s})\, ds\\
&&\hskip4cm
=\frac{2}{\te^2}\,a_t,
\end{eqnarray*}
and the statement in Proposition \ref{!!!} can be formulated as follows
\begin{proposition}
\label{occu}
{\sl Let $Z,\ Z_0=x,$ and  $X,\ X_0=1/(c\,{\rm e}^{\,\te x}+1),$
be as above. Then
\begin{equation}
\label{o1}
\frac {\te^2}{2}
\int_0^\infty\frac{c\,\exp(\te Z_s)}{(c\,\exp(\te Z_s)+1)^2}\, ds
=
H_0(X)\quad {\rm a.s.}
\end{equation}
}
\end{proposition}
%\quad{\mathop=^{\rm{(d)}}}\quad
\hskip1cm
%f(X_t)-f(X_0)=\int_0^t f'(X_s)\,dB_s + \int_0^t f'(X_s)(a-b\,X_s)\, ds
%+\frac{1}{2}\int_0^tf''(X_s)\, 2X_s(1-X_s)\, ds

As explained in Section \ref{sec1}, we can use the Feynman--Kac
method to deduce that the Laplace-transforms of the functionals in
(\ref{o1}) are equal. However, because the present case is not covered
by Proposition \ref{rk!}, we formulate here 
a result connecting
the solutions of the
hypergeometric differential
equation to the solutions of the equation
induced by the generator of $Z$ with the
potential term as on the left hand side of (\ref{o1}).
%This result is needed to find the explicit form of the Laplace
%transform of the distributions of the random variables in (\ref{o1}).

\begin{proposition}
\label{prop1}{\sl 
Let $z(x)=z(\al,\be,\ga,x)$, $x\in (0,1)$, be an 
arbitrary solution of the hypergeometric differential equation
\begin{equation}
\label{11}
x(1-x)z''(x)+(\ga-(\al+\be+1)x)z'(x)-\al\be z(x)=0.
\end{equation}
Then for $c\ge 0$ and $\te\in\R$
the function $q(x)=z\big(y(x)\big)$
with 
$$
y(x):= (c\,{\rm e}^{\,\te x}+1)^{-1}
$$
satisfies for $x\in (-\infty,\infty)$ the equation
\begin{equation}
\label{121}
q''(x)+\te\Big(1-\ga+\frac{\al+\be-1}
{c\,{\rm e}^{\,\te x}+1}\Big)q'(x)
-\frac{\te^2\al\be c\,{\rm e}^{\,\te x}}{\big(c\,{\rm e}^{\,\te x}+1\big)^2}q(x)=0.
\end{equation}
}
\end{proposition}

\begin{proof} It is sufficient to prove this statement
for $a=1$ since one can shift the
argument. Differentiating the composition of the functions we have (the
corresponding arguments are omitted)
$$
q'=-\frac{\te{\rm e}^{\,\te x}}{\big({\rm e}^{\,\te x}+1\big)^2}z',\qquad
q''=\frac{\te^2 {\rm e}^{\,2\te x}}{\big({\rm e}^{\,\te x}+1\big)^4}z''
+\frac{\te^2 {\rm e}^{\,\te x}\big({\rm e}^{\,\te x}-1\big)}
{\big({\rm e}^{\,\te x}+1\big)^3}z'.
$$
Taking in (\ref{11}) the argument $y(x)$
instead of $x$, one obtains
$$\frac{{\rm e}^{\,\te x}}{\big({\rm e}^{\,\te x}+1\big)^2}z''
+\Big(\ga-\frac{\al+\be+1}{{\rm e}^{\,\te x}+1}\Big)z'-\al\be z=0. $$
Finally, from these relations we have
$$
q''(x)=\frac{\te^2{\rm e}^{\,\te x}}{\big({\rm e}^{\,\te x}+1\big)^2}
\Big(\Big(\frac{\al+\be+1}{{\rm e}^{\,\te x}+1}-\ga
+\frac{{\rm e}^{\,\te x}-1}{{\rm e}^{\,\te x}+1}\Big)z'+\al\be z\Big)$$
$$=-\te\Big(\frac{{\rm e}^{\,\te x}+\al+\be}
{{\rm e}^{\,\te x}+1}-\ga\Big)q'(x)+
\frac{\te^2\al\be {\rm e}^{\,\te x}}{\big({\rm e}^{\,\te x}+1\big)^2}q(x) $$
proving the equation (\ref{121}).
\end{proof}

\begin{remark}
\label{converse}
{\rm Notice also the
 converse of Proposition \ref{prop1}: if $q$ is a solution of (\ref{121}) then
$z(x)=q(y^{-1}(x)),$ 
where 
$$
y^{-1}(x):=\frac{1}{\te} \log\Big(\frac{1-x}{cx}\Big),\quad 0<x<1,
$$
is a solution of (\ref{11}) (cf. (\ref{Z})).
}\end {remark}

\subsection{Functionals of Jacobi diffusion; Case 2}
\label{sec23}

Assume that  $X$ is a Jacobi diffusion with $\ga=\al+\be\geq 2$
and recall that in this case $X$ is recurrent.
By Ito's formula %for $t<H_0(X)$
\begin{eqnarray}
\label{e11}
\nonumber
&&
\log\left(1- X_t\right) -\log\left(1- X_0\right)
\\
\nonumber
&&\hskip.5cm
=
-\sqrt{2}\int_0^t \sqrt{\frac{X_s}{1-X_s}}\,dB_s\\
&&\hskip2cm
\nonumber
- \int_0^t\frac{\ga-(\al+\be+1)X_s}{1-X_s} \, ds
-\frac{1}{2}\int_0^t\frac{2\,X_s(1-X_s)}{(1-X_s)^2} \, ds\\
\nonumber
&&
\\
&&\hskip.5cm
=-\sqrt{2}\int_0^t \sqrt{\frac{X_s}{1-X_s}}\,dB_s
-\int_0^t\frac{X_s}{1-X_s}\left( \frac{\ga}{X_s}-\al-\be\right)
\, ds
\\
\nonumber
&&
\\
\nonumber
&&\hskip.5cm
=-\sqrt{2}\int_0^t \sqrt{\frac{X_s}{1-X_s}}\,dB_s
- \int_0^t\frac{X_s}{1-X_s}\frac{(\al+\be)(1-X_s)}{X_s}
\, ds.
\end{eqnarray}
Let
$$
\widehat A_t:=\frac{2}{\te^2}\int_0^t \frac{X_s}{1-X_s}\,ds,
$$
and define its inverse $\widehat a_t:=\inf\{s\,:\, \widehat A_s>t\}.$
By continuity, $\widehat A_t< \infty$ and $\widehat a_t<\infty$
for all $t$ and $\widehat A_t\to \infty$ and $\widehat a_t\to \infty$
as $t\to\infty;$ hence, we are in the case covered by Proposition \ref{!!}.1.
Instead of simply refering to Proposition \ref{!!}.1, we give some
details. Firstly
\begin{equation}
\label{alpha11}
\widehat a_t=\frac{\te^2}{2}\int_0^t \frac {1-X_{a_s}}{X_{a_s}}\,ds,
\end{equation}
and the process given by
\begin{equation}
\label{hatZ}
\widehat Z_t:=-\frac{1}{\te}\log (1-X_{\widehat a_t})
\end{equation}
is well defined for all $t\geq 0.$
From (\ref{e11}) we obtain
\begin{eqnarray*}
&&\widehat Z_t-\widehat Z_0
= \widehat B^\circ_t+
\frac {\te}{ 2}\int_0^t
\frac{\al+\be}{\exp(\te\,\widehat Z_s)-1}\, ds
\end{eqnarray*}
where $\widehat B^\circ $ is a Brownian motion. 
From (\ref{alpha11}) and (\ref{hatZ}) we obtain
%$\widehat Z_t\to +\infty$ a.s. and from (\ref{alpha11})
\begin{eqnarray*}
&&
\widehat a_t=\frac{\te^2}{2}\int_0^t \frac{1-X_{a_s}}{X_{a_s}}\, ds%\\
%&&\hskip.55cm
=\frac{\te^2}{2}\int_0^t\left(\exp(\te\,\widehat Z_s)-1\right)^{-1}ds
\end{eqnarray*}
For $y>0$ introduce
$$
H_y(\widehat Z):=\inf\{t\,:\, \widehat Z_t=y\}.
$$
Now, from (\ref{hatZ}),
$$
y=\widehat Z_{H_y(\widehat Z)}=-\frac{1}{\te}\log \left(1-X_{\widehat a_{H_y(\widehat Z)}}\right)
$$
and, hence,
$$
\widehat a_{H_y(\widehat Z)}=\inf\{t\,:\, X_t=1-{\rm e}^{-\te\,y}\}.
$$
Consequently, we arrive to the result (cf. Proposition \ref{!!}.1)

\begin{proposition}
\label{occu2}
{\sl Let $X$ be a Jacobi diffusion 
with $\ga=\al+\be\geq 2$ and $\widehat Z$ as defined in (\ref{hatZ})
with $\widehat Z_0=x.$ Then for $y>0$
\begin{equation}
\label{o11}
\frac{\te^2}{2}\int_0^{H_y(\widehat Z)}
\left(\exp(\te\,\widehat Z_s)-1\right)^{-1}ds
=
H_{y^*}(X).
\end{equation}
where  $y^*=1-{\rm e}^{-\te\,y}$ and $X_0=1-{\rm e}^{-\te\,x}.$
}
\end{proposition}

It is useful (again) to give the corresponding result 
for differential equations.
Notice that in this result we do not ìnterprete solutions
probabilistically and, therefore, can formulate a more general
statement without any restrictions on the values of the parameters.

\begin{proposition}
\label{prop3}{\sl
Let $z(x)=z(\al,\be,\ga,x)$, $x\in (0,1)$, be an
arbitrary solution of the hypergeometric differential equation (\ref{11}).
Then for $\te>0$ the function $s(x)=z\left(1-{\rm e}^{-\te x}\right)$
satisfies for $x\in (0,\infty)$ the equation
$$
s''(x)+\te\Big(\frac{\ga\,{\rm e}^{\,\te x}}{{\rm e}^{\,\te x}-1}
-\al-\be\Big)s'(x)
-\frac{\te^2\al\be }{{\rm e}^{\,\te x}-1}s(x)=0.
$$
}
\end{proposition}

\begin{proof} The proof is analogous to the proof of Proposition \ref{prop1}.
We have 
$$
s'=\te {\rm e}^{\,-\te x}z',\qquad
s''=\te^2 {\rm e}^{\,-2\te x}z''-\te^2 {\rm e}^{\,-\te x} z'.
$$
Taking in (\ref{gauss}) the argument $1-{\rm e}^{\,-\te x}$
instead of $x$ yields
$$ {\rm e}^{\,-\te x}z''
-\Big(\al+\be+1-\frac{\ga}{1-{\rm e}^{\,-\te x}}\Big)z'
-\frac{\al\be}{1-{\rm e}^{\,-\te x}}z=0. $$
Consequently,
\begin{eqnarray*}
s''(x)&=&\te^2{\rm e}^{\,-\te x}
\Big(\al+\be+1-\frac{\ga}{1-{\rm e}^{\,-\te x}}-1\Big)z'
+\frac{\te^2\al\be {\rm e}^{\,-\te x}}{1-{\rm e}^{\,-\te x}}z
\\
&=&\te\Big(\al+\be-\frac{\ga}{1-{\rm e}^{\,-\te x}}\Big)s'(x)+
\frac{\te^2\al\be {\rm e}^{-\te x}}{1-{\rm e}^{\,-\te x}}s(x),
\end{eqnarray*}
proving the claim.
\end{proof}

\begin{remark}
\label{1-x2}
{\rm {\bf 1.} We have also
\begin{eqnarray*}
&&
\log\left(X_t\right) -\log\left(X_0\right)
\\
&&\hskip1cm
=\sqrt{2}\int_0^t \sqrt{\frac{1-X_s}{X_s}}\,dB_s
+\int_0^t\frac{1-X_s}{X_s}\left(\al+\be- 
\frac{\al+\be+1-\ga}{1-X_s}\right)
\, ds.
\end{eqnarray*}
Choosing $\ga=1$ and letting  
$$
A^\circ_t:=\frac{2}{\te^2}\int_0^t \frac{1-X_s}{X_s}\,ds,
$$
it is seen that the process
$$
Z^\circ_t:=-\frac{1}{\te}\log \left(X_{a^\circ_t}\right),
$$
where $a^\circ$ is the inverse of $A^\circ$, satisfies 
the same SDE as $\widehat Z.$

\noindent
{\bf 2.}
Let $z(x)=z(\al,\be,\ga,x)$, $x\in (0,1)$, be an
arbitrary solution of the hypergeometric differential equation
(\ref{11}).
%\begin{equation}
%\label{11}
%x(1-x)z''(x)+(\ga-(\al+\be+1)x)z'(x)-\al\be z(x)=0.
%\end{equation}
Then for $\te\ge 0$ %and $\te\in\R$
the function $q(x)=z\big({\rm e}^{-\te x}\big)$
satisfies for $x>0$ the equation
\begin{equation}
\label{12}
\,q''(x)+\te\Big(1-\ga+\frac{\al+\be+1-\ga}
{{\rm e}^{\,\te x}-1}\Big)q'(x)
-\frac{\te^2\al\be }{{\rm e}^{\,\te x}-1}\,q(x)=0.
\end{equation}
} 
\end{remark}

\section{Perpetual integral functional of BM($\mu$)}
\label{sec3}

Our first main result completes, in a sense, the result in \cite{salminenyor03b}, see
also \cite{salminenyor04}, concerning the translated Dufresne's
functional but, moreover, it gives Laplace transforms for many new
perpetual integral functionals. The functional we analyze is
$$
I_\infty(p,q):=\int_0^\infty
\Big(\,
\frac{p}{c\,\exp(\te\,B^{(\mu)}_s)+1}+
\frac{q}{\big(c\,\exp(\te\, B^{(\mu)}_s)+1\big)^{2}}\Big)
\, ds.
$$
The notation $I_t(p,q)$ is used when the integration is from 0 to $t$ and
$I_{H_y}(p,q)$ when $t$ equals
$$
H_y:=\inf\{s\,:\, B^{(\mu)}_s=y\},
$$
the first hitting time of $y.$
\begin{theorem}
\label{main}
{\sl
Let $B^{(\mu)}$ be a Brownian motion with drift $\mu>0$ started from
$x.$ Then for $c>0,\,\te>0,\,  p\geq 0$ and $p+q\geq 0$
\begin{eqnarray}
\label{a101}
&&
\nonumber
\E_x\Big(\exp\Big(-I_\infty(p,q)\Big)\Big)
\\
&&\hskip2cm
= K
\ v(x)^k\ F(\alpha,\beta,\alpha+\beta+2\mu/\te\,;\,v(x))
%\\
%&&\hskip2cm
%\nonumber
%= K\ v(x)^k\
%G(\alpha,\beta,1-2\mu/\te\,;\,1-v(x)),
%
\end{eqnarray}
where $F$ is 
%and $G$ are 
Gauss' hypergeometric function as
defined in (\ref{F}), %and (\ref{G}), respectively, 
and
$$
k=(\alpha+\beta-1)/2,\qquad v(x)=\frac{c\,\exp(\te\,x)}{c\,\exp(\te\,x)+1}
$$
\begin{equation}
\label{K}
K=\frac{\Gamma(\alpha+2\mu/\te)\Gamma(\beta+2\mu/\te)}
{\Gamma(\alpha+\beta+2\mu/\te)\,\Gamma(2\mu/\te)}
\end{equation}
and
\begin{equation}
\label{alpha}
\al=\frac12-\mu/\te+\sqrt{\mu^2+2(p+q)}/\te+\frac12\sqrt{1+8q/\te^2},
\end{equation}
\begin{equation}
\label{beta}
\be=\frac12-\mu/\te+\sqrt{\mu^2+2(p+q)}/\te-\frac12\sqrt{1+8q/\te^2}.
\end{equation}
}
\end{theorem}
%\begin{remark}
%\label{alt}
%Using the function $G$ defined in (\ref{G}) the function $F$ on the
%right hand side of (\ref{a101}) can be expressed as
%$$
%F(\alpha,\beta,\alpha+\beta+2\mu/\te\,;\,v(x))
%= G(\alpha,\beta,1-2\mu/\te\,;\,1-v(x)).
%$$
%\end{remark}
\begin{proof}
Recall, e.g., from \cite{salminenyor04} that
\begin{eqnarray*}
&&\Psi(x):=\E_x\Big(\exp\Big(-I_\infty(p,q)
\Big)\Big)
\\
&&\hskip1.2cm
= \E_x\Big(\exp\Big(-\int_0^\infty
\Big(\,
p\,\Big(1-v(B^{(\mu)}_t)\Big) +
q\,\Big(1-v(B^{(\mu)}_t)\Big)^2\Big)
\, dt\Big)\Big)
\end{eqnarray*}
is the unique positive bounded function such that
\begin{equation}
\label{ode}
\frac{1}{2}\Psi^{\, ''}(x)+\mu\,\Psi^{\, '}(x)
-\Big(p\,(1-v(x))+q\,(1-v(x))^2\Big)\,\Psi(x)
=0
\end{equation}
and $\lim_{x\to+\infty}\Psi(x)=1.$ We remark that in
\cite{salminenyor04} it is required $\Psi$
to be increasing, but from the proof therein it is clear
that an equivalent requirement is boundedness.
To find $\Psi$ we use Girsanov's theorem and Proposition \ref{occu}.
Firstly, recall that the process $Z$ as defined in (\ref{Z})
solves the SDE
$$
dZ_t=dB^\circ_t+\frac{\te}{2}\Big(1-\ga+(\al+\be-1)(1-v(Z_t))
\Big)\,dt.
$$
Choose now $\ga$ such that $\mu=\te(1-\ga)/2.$
By Girsanov's theorem,
the measures $\P^Z$ and $\P^{(\mu)}$  induced  by $Z$ and $B^{(\mu)},$
respectively, and defined in the space of
  continuous functions $\omega:\R_+\mapsto\R$ are locally absolutely
  continuous with the Radon--Nikodym derivative given by
\begin{equation}
\label{abs}
{\P^Z}\, |_{\cF_t} =
\exp\left(D_t\right)\, {\P^{(\mu)}}\, |_{\cF_t},
\end{equation}
where $\cF_t$ is the $\sigma$-algebra generated by the coordinate
process up to time $t$ and 
\begin{equation}
\label{RN}
D_t:=\int_0^t f(\omega_s)d\omega_s-\frac 12 \int_0^t\Big(
f^2(\omega_s)+ 2\mu\, f(\omega_s)\Big)ds
\end{equation}
with
$$
f(\omega_s)=\te\,k\,(1-v(\omega_s)).%\left(\al+\be-1\right).
$$
Observe
$$
V(x):=\int^x(1-v(y))\,dy=\frac1\te\int^x\frac {v'(y)}{v(y)}\,dy =
\frac1\te\log v(x).
$$
For the first integral in (\ref{RN}) we obtain by Ito's formula
\begin{eqnarray*}
&&\int_0^t(1-v(\omega_s))\,d\omega_s= V(\omega_t)-V(\omega_0)+\frac
12\int_0^t v'(\omega_s)\,ds\hskip2cm\\
&&\hskip3.2cm
= V(\omega_t)-V(\omega_0)+\frac
\te2\int_0^t v(\omega_s)(1-v(\omega_s))\,ds\\
&&\hskip3.2cm
= V(\omega_t)-V(\omega_0)-\frac
\te2\int_0^t (1-v(\omega_s))^2\,ds\\
&&\hskip7cm +\frac
\te2\int_0^t (1-v(\omega_s))\,ds.
\end{eqnarray*}
Letting (cf. Proposition \ref{occu})
$$
I_t:= \frac{\te^2}{2} 
\int_0^t\frac{c\,\exp(\te \omega_s)}{(c\,\exp(\te \omega_s)+1)^2}\, ds
$$
we have 
\begin{eqnarray*}
&&I_t
=\frac{\te^2}{2} \int_0^t v(\omega_s)(1-v(\omega_s))\, ds\\
&&\hskip.5cm
=\frac{\te^2}{2} \int_0^t (1-v(\omega_s))\, ds
-\frac{\te^2}{2} \int_0^t (1-v(\omega_s))^2\, ds.
\end{eqnarray*}
The absolute continuity relation (\ref{abs}) yields
\begin{equation}
\label{21}
\E^Z_x\Big(\exp(-r\, I_t)\Big)=
\E^{(\mu)}_x\Big(\exp(D_t-r\, I_t)\Big).
\end{equation}
Straightforward but lengthy computations show that choosing $r=\al\be$
the claimed expressions (\ref{alpha}) and (\ref{beta}) for $\al$ and $\be,$ respectively, are such that
\begin{eqnarray}
\label{22}
&&
\nonumber
\E^{(\mu)}_x\Big(\exp(D_t-r\, I_t)\Big)\\
&&\hskip1cm
=\E_x\Big[\exp\Big(\te k \,V(B^{(\mu)}_t)-
\te k\, V(x)\Big)
\exp\Big(-I_t(p,q)\Big)
\Big].
\end{eqnarray}
We remark that $r=\al\be$ can attain both positive and negative
values. The equalities (\ref{21}) and (\ref{22})
hold also for the first hitting time $H_y$ with $y>x,$
i.e.,
\begin{eqnarray*}
&&
\E^Z_x\Big(\exp(-r\, I_{H_y})\Big)=
\left(\frac{v(y)}{v(x)}\right)^k
\E_x^{(\mu)}\Big[\exp\Big(-I_{H_y}(p,q)\Big)\Big].
\end{eqnarray*}
Letting $y\to +\infty$ (remember that $Z_t\to +\infty$ and  $B^{(\mu)}_t\to +\infty$
as $t\to +\infty$) we obtain by monotone convergence and Proposition \ref{occu}
\begin{equation}
\label{13}
\Psi(x)%\E_x^{(\mu)}\Big[\exp\Big(-I_{\infty}(p,q)\Big)\Big]
= (v(x))^k \E^X_{x^*}\Big[\exp\Big(-\al\,\be\,H_0(X)\Big)\Big],
\end{equation}
where $x^*:=1-v(x).$ Recall that $\psi(x)\to 1$ as $x\to +\infty.$
This can also be verified from (\ref{13})
by observing that as $x\to +\infty$ then $x^*\to 0$ and
the right hand side of (\ref{13})
tends to 1 because $0$ is an exit
(or killing) boundary point for $X.$
Notice also that if $p+q=0$ (see Example 5.3) (in this case $\al\be>0$)
then $k=0$ and we have
$$
0<\lim_{x\to
  -\infty}\Psi(x)
=
\lim_{x^*\to 1}\E^X_{x^*}\Big[\exp\Big(-\al\,\be\,H_0(X)\Big)\Big].
$$
Also if $p+q>0$ then $\Psi(x)\to 0$ as $x\to -\infty.$ In
this case, because $v(x)\to 0$ we cannot without further analysis
claim that
$$
\lim_{x^*\to 1}\E^X_{x^*}\Big[\exp\Big(-\al\,\be\,H_0(X)\Big)\Big]<+\infty
$$
in case $\al\be<0$ (but this will follow from our discussion !). 
For the general theory of diffusions we know that the function
$$
\psi^X_{\al\be}(x^*):=\E^X_{x^*}\Big[\exp\Big(-\al\,\be\,H_0(X)\Big)\Big]
$$
is a solution of the hypergeometric differential equation

(Clearly, $\psi^X_{\al\be}$ as a function of $x^*$ is decreasing 
if $\al\be>0,$ and increasing 
if $\al\be<0.$ Hence, because $1-v(x)$ is decreasing,
$\psi^X_{\al\be}$ as a function of $x$ is increasing 
if $\al\be>0,$ and decreasing 
if $\al\be<0.$ In the latter case, multiplying $\psi^X_{\al\be}$
with the increasing function $v(x)^k$ makes the product increasing.)

Let $F$ be the hypergeometric function as defined in (\ref{F})
with  $\al$ and $\be$ as in (\ref{alpha}) and (\ref{beta}),
respectively, and $\ga=1-(2\mu/\te).$ Notice that 
$\al$ and $\be$ can be conjugate complex numbers and
$$
\al+\be>1,\quad{\rm and}\quad  \ga<1.
$$ 
Moreover, as is easily seen,
$$
{\rm Re}(\be+1-\ga)={\rm Re}(\be+(2\mu/\te))>0.
$$
From Section \ref{sec1} it now follows (cf. (\ref{G2})) that the function
\begin{eqnarray*}
&&
%G(\alpha,\beta,\gamma;x):=
x\mapsto F(\alpha,\beta,\al+\be+ 1 -\gamma;1-x)%\\
%&&\hskip2.4cm
=F(\alpha,\beta,\al+\be+ (2\mu/\te);1-x)
\end{eqnarray*}
is a bounded non-negative solution of the hypergeometric differential
equation. Proposition \ref{prop1} and formula (\ref{13}) yield
$$
\Psi^\circ(x):= (v(x))^k\, 
F(\alpha,\beta,\al+\be+ (2\mu/\te);v(x))
%G(\alpha,\beta,\gamma;1-v(x))
$$
is a nonnegative bounded solution of (\ref{ode}).
Consequently, by the
uniqueness we have
$$
\Psi(x)=\Psi^\circ(x)/\Psi^\circ(+\infty),
$$
where, from (\ref{Glimit})
\begin{eqnarray*}
&&\Psi^\circ(+\infty):=\lim_{x\to +\infty}\Psi^\circ(x)\\
&&\hskip1.8cm
=F(\alpha,\beta,\al+\be+ (2\mu/\te);1)\\
%=G(\alpha,\beta,\gamma;0)\\
&&\hskip1.8cm
=\frac{\Gamma(\alpha+\beta+(2\mu/\te))\,\Gamma(2\mu/\te)}
{\Gamma(\al+(2\mu/\te))\,\Gamma(\beta+(2\mu/\te))},
\end{eqnarray*}
as claimed.
\end{proof}

%h

\begin{example}
\label{sal_yor}
{\rm Choosing in (\ref{a101}) $p=0,$ $q=\gamma/a^2,$ $\theta=1,$ and $c=1/a$
gives 
\begin{eqnarray*}
&&\E_x\Bigr(\exp\Bigr(-\gamma\int_0^\infty
{{(a+\exp(B^{(\mu)}_s))^{-2}}}\, ds
\Bigr)\Bigr)\\
&&\hskip2cm
=
K\
v(x)^{\,(\sqrt{a^2\mu^2+2\gamma}-a\mu)/a}\
F(\alpha,\beta,\al+\beta +2\mu\,;v(x)),
\end{eqnarray*}
where 
$$
\alpha=\frac 12 -\mu +\sqrt{\mu^2+\frac{2\gamma}{a^2}}+
\sqrt{\frac 14+\frac{2\gamma}{a^2}},\quad
\beta=\frac 12 -\mu +\sqrt{\mu^2+\frac{2\gamma}{a^2}}-
\sqrt{\frac 14+\frac{2\gamma}{a^2}},
$$
$K$ is given by (\ref{K}) with $\alpha$ and $\beta$ as above, and    
$$
v(x)=
{\exp(x)\over{a+\exp(x)}}.
$$
Notice that if $\mu=1/2$ then $\beta=0$ and because
$F(\alpha,0,\gamma;x)=1$ for all $|x|<1$ (see (\ref{gauss}))
we obtain the result in \cite{salminenyor04} (see also \cite{salminenyor03b}):
\begin{eqnarray*}
%\label{01}
%\nonumber
&&\E_x\Bigr(\exp\Bigr(-\gamma\int_0^\infty
{{(a+\exp(B^{(1/2)}_s))^{-2}}}\, ds
\Bigr)\Bigr)
%&&\hskip3cm
=v(x)^{(2\,a)^{-1}\,(\sqrt{a^2+8\gamma}-a)}.
\end{eqnarray*}
%\Bigl({\exp(x)\over{a+\exp(x)}}
%\Bigr)
From Proposition \ref{!!}.2 it follows that 
$$
\int_0^\infty
{{(a+\exp(B^{(\mu)}_s))^{-2}}}\, ds
=H_0(Z)\quad {\rm a.s.}
$$
where $Z$ satisfies the SDE
$$
dZ_t=dB_t+\left((\mu-\frac 12)\frac{a\,{\rm e}^{aZ_t}}
{1-{\rm e}^{aZ_t}}+a\mu\right)dt.
$$ 
The function $g$ is in this case
$$
g(x):=-\frac 1a\,\ln(1+a{\rm e}^{-x}).
$$
}
\end{example}
%here
\begin{example}
{\rm 
Let in (\ref{a101}) $p+q=0$, $p=4\rho\ge 0$,
$\te=2$, and $c=1.$  Then  
$$
\alpha=\frac 12 +\frac 12\sqrt{1-8\rho},\quad
\beta=\frac 12 -\frac 12\sqrt{1-8\rho},\quad k=0,
$$
$$
K=
\frac{\Gamma\left(\mu+\al\right)\ \Gamma\left(\mu+\beta\right)}
{\Gamma(\mu)\ \Gamma(\mu+1)},
$$
and we obtain the result (see Vagurina \cite{vagurina04})
\begin{eqnarray*}
&&\E_x\left(\exp\Bigl(-\rho\,\int_0^\infty
 \ch^{-2}(B^{(\mu)}_s)\ ds\Bigr)\right)
=K\,F(\al,\beta,1+\mu;v(x)).
%&&\hskip1.5cm
%\times
%G\Big((1+\frac12\sqrt{1-8\rho})/2,
%(1-\sqrt{1-8\rho})/2,1-\mu;(e^{2x}+1)^{-1} \Big).
\end{eqnarray*}
By Proposition \ref{!!}.2 
\begin{equation}
\label{d1}
\int_0^\infty \frac{ds}{\ch^2(B^{{(}\mu{)}}_s)}=H_\pi(Z), \qquad a.s,
\end{equation}
with $Z$ determined via the SDE
$$
dZ_t=dB_t+\left(\frac12\,\ctg Z_t+\frac{\mu}{\sin Z_t}\right)\,dt,\qquad
Z_0=2\arctg\exp(B^{{(}\mu{)}}_0). 
$$
In this example 
$$ 
g(x):=2\arctg {\rm e}^x,\quad g'(x)=\frac1{\cosh x},
$$ 
and, therefore,  
$$
G(g^{-1}(x))=\frac12\, \ctg x+\frac{\mu}{\sin x}
=\frac12\,\left(\mu-\frac12\right)\,\tg\frac x2+
\frac12\left(\mu+\frac12\right)\ctg\frac x2. 
$$
Further,
$$
\E_z(\exp(-\rho H_\pi(Z)))=
K\,F(\al,\beta,1+\mu;\sin^2(z/2))
$$
}
\end{example}

\begin{example}
{\rm Take $p=hc$ and  $q=a^2\,c^2/2$ in the defintion of $I_\infty(p,q).$
Then, letting $c\to +\infty,$ we obtain by monotone convergence (using 
(\ref{e01}))
\begin{eqnarray*}
&&
\lim_{c\to +\infty}
\int_0^\infty
\Big(\,
\frac{hc}{c\,\exp(\te\,B^{(\mu)}_s)+1}+
\frac{a^2}{2}\ \frac{c^2}{(c\,\exp(\te\, B^{(\mu)}_s)+1\big)^{2}}\Big)
\, ds\\
&&\hskip3cm
=
\int_0^\infty
 \left(h\, \exp(-\te B^{(\mu)}_s)+
 \frac{a^2}{2}\,\exp(-2\te B^{(\mu)}_s)\right)ds
\end{eqnarray*}
The Laplace transform of the functional on the right hand side 
is given in \cite{borodinsalminen02}
formula 2.1.30.3(2) p. 292 and can written as
\begin{eqnarray}
\label{cf}
\nonumber
&&
\E_x\left(\exp\left(-\int_0^\infty
 \left(h\,\exp(-\te B^{(\mu)}_s)+
 \frac{a^2}{2}\,\exp(-2\te B^{(\mu)}_s)\right)ds\right)\right)
\\
&&\hskip2cm
=
\frac{\Gamma(1/2+\mu/\te+h/a\te)}{\Ga(2\mu/\te)}\exp\Big(-\frac c\te\  {\rm e}^{-\te x}\Big)
\\
\nonumber
&&\hskip3.5cm
\times\
U\Big(\frac12-\frac{\mu}{\te}+\frac{h}{a\te},1-\frac{2\mu}{\te},
\frac{2a}{\te}\ {\rm e}^{-\te x}\Big),
\end{eqnarray}
where the Kummer function $U$ is connected to the Whittaker function $W$ via
$$
W_{n,m}(x)=W_{n,-m}(x)=x^{-m+1/2}\ {\rm e}^{-x/2}\,U(-m-n+1/2,1-2m,x),
$$
see Abramowitz and Stegun \cite{abramowitzstegun70}.
To obtain the formula (\ref{cf}) from (\ref{a101})
observe first that
as $c\to +\infty$ we have
$\al\simeq 2ac/\te,$
$$
\be\to \frac12-\frac\mu\te+\frac{h}{a\te},
$$
and
$$
\Big(1-\frac1{c\,{\rm e}^{\,\te x}+1}\Big)^{ac/\te} \to\exp\Big(-\frac a\te\, {\rm e}^{-\te x}\Big).
$$
and recall that for Re$\ga<1$ and Re$(\al+\be-\ga)>-1$
$$
\lim\limits_{\al\to \infty}\frac{\Ga(\al+1-\ga)}{\Ga(\al+\be+1-\ga)}
F\bigl(\al,\be,\al+\beta+1-\ga;1-\frac x\al\bigr)
=U(\be,\ga,x),
$$
see Abramowitz and Stegun \cite{abramowitzstegun70}.
}
\end{example}

\section{Perpetual integral functional of BES(3)}
\label{sec4}

Let $R=\{R_t\,:\, t\geq 0\}$ denote a 3-dimensional Bessel process (or, equivalently, of index $1/2$). The generator of $R$ is
$$
\cG^Ru={1\over 2}{d^2u\over{dx^2}}
+{1\over {x}}\,{du\over{dx}}.
$$
In the next theorem we give the Laplace transform of the functional.
$$
\widehat I_\infty(p,q):=\int_0^\infty
\Big(\,
\frac{p}{\exp(\te\,R_s)-1}+
\frac{q}{(\exp(\te\, R_s)-1\big)^{2}}\Big)
\, ds.
$$

\begin{theorem}
\label{main2}
{\sl
Let $R$ be a 3-dimensional Bessel process started from
$x>0.$ Then for $\te>0,\,  p\geq 0$ and $q\geq 0$
\begin{eqnarray}
\label{b1}
\nonumber
&&\E_x\Big(\exp\Big(-\widehat I_\infty(p,q)\Big)\Big)\\
&&
\hskip2cm
=\frac{(1-{\rm e}^{-\te x})^{\ga/2}\, \Ga(\al)\,\Ga(\be)}{x\,\te\,\Ga(\ga)}\
F\big(\al,\be,\ga,1-{\rm e}^{-\te x}\big)
\\
\nonumber
&&\hskip2cm
=\frac{(1-{\rm e}^{-\te x})^{\ga/2}}{x\,\te}\,
\int_0^1 t^{\al-1}(1-t)^{\be-1}(1-tx)^{\al-\ga}dt,
\end{eqnarray}
where 
\begin{equation}
\label{aa1}
\al=\frac12+\frac12\sqrt{1+8q/\te^2}+\frac1\te\sqrt{2(q-p)}
\end{equation}
\begin{equation}
\label{bb1}
\be=\frac12+\frac12\sqrt{1+8q/\te^2}-\frac1\te\sqrt{2(q-p)},
\end{equation}
and $\ga=\al+\be=1+\sqrt{1+8q/\te^2}.$
}
\end{theorem}

\begin{proof}
To start with, let  $\widehat Z^\circ$ be a diffusion with the generator
$$
\cG^\circ u={1\over 2}{d^2u\over{dx^2}}
+h(x)
\,{du\over{dx}},
$$
where $\ga>\al+\be\geq 0$ and 
$$
h(x):=\frac \te 2\left(\frac{\ga\,{\rm e}^{\,\te\,x}}
{{\rm e}^{\,\te\,x}-1}-(\al+\be)\right).
$$
Notice that if $\ga=\al+\be$ then
$\widehat Z^\circ$ coincides with the diffusion $\widehat Z$
introduced in Section \ref{sec23}. The measures induced by $\widehat Z^\circ$ and
by $R$ in the canonical space of continuous functions are absolutely continuous
with respect to each other when restricted to the $\sigma$-algebra
$\cF_t$
 generated by the co-ordinate mappings up to a fixed but 
arbitrary time $t.$ Let $\P^{\widehat Z}_x$ and  $\P^{R}_x$ be 
the measures associated with $\widehat Z^\circ$ and $R,$ 
respectively, when both processes are started from $x>0.$ Then
\begin{equation}
\label{abs2}
\P^{\widehat Z}_x\, |_{\cF_t} =
%{\rm e}^{\widehat D_t}\, \P^{R}_x\, |_{\cF_t},
\exp\big(\widehat D_t\big)\, \P^{R}_x\, |_{\cF_t},
\end{equation}
where the exponent $\widehat D_t$ in the Radon--Nikodym derivative is given by
\begin{equation}
\label{RN2}
\widehat D_t=-\int_0^t(\omega_s^{-1}
-h(\omega_s))
\, d\omega_s
+\frac 12\int_0^t(\omega_s^{-2}-h^2(\omega_s))
\, ds
\end{equation}
Consider the stochastic integral term in (\ref{RN2}). Under the
measure $\P_x^R$ the co-ordinate process is the Bessel process $R$
started from $x.$
Recalling that the quadratic variation of $R$ is $t,$ we obtain by
Ito's formula 
\begin{enumerate}
\item \quad  ${\displaystyle -\int_0^t R_s^{-1}
\, dR_s =\log\left(\frac{x}{R_t}\right) - 
\frac 12\int_0^tR_s^{-2}\,ds}$
\item\quad $
{\displaystyle 
\int_0^t h(R_s)\,dR_s
=\frac\ga 2\log\left(\frac{{\rm e}^{\,\te\,R_t}-1}{{\rm
  e}^{\,\te\,x}-1}\right)
-\frac{\te\,(\al+\be)} 2(R_t-x)}
$
$$
\hskip3cm {\displaystyle +\frac {\te^2\,\ga} {4}
\int_0^t\frac{\,{\rm e}^{\,\te\,R_s}}
{({\rm e}^{\,\te\,R_s}-1)^2}\,ds.
}
$$
\end{enumerate}
Next consider the bounded variation part in (\ref{RN2}). We have 
\hfill\break\hfill

\noindent
\hskip.2cm
(iii)\hskip.6cm
${\displaystyle -\frac 12\int_0^t h^2(R_s)\, ds= 
-\frac {\te^2\,\ga^2} {8}
\int_0^t\frac{\,{\rm e}^{\,2\,\te\,R_s}}
{({\rm e}^{\,\te\,R_s}-1)^2}\,ds
- \frac {\te^2\,(\al+\be)^2} {8}\, t}\hskip2cm$
\hfill\break\hfill

 $\hskip6cm{\displaystyle 
+\frac {\te^2\,\ga\,(\al+\be)} {4}\, 
\int_0^t\frac{\,{\rm e}^{\,\te\,R_s}}
{{\rm e}^{\,\te\,R_s}-1}\,ds}.$

\vskip.5cm
\noindent
Because $R_t\to +\infty$  as $t\to +\infty,$ it is seen from (ii) and
(iii) above that taking $\ga=\al+\be$ will lead to a remarkably
simple special case. Indeed, if $\ga=\al+\be$

\begin{eqnarray*}
&&\widehat D_t=  \log\left(\frac{x}{R_t}\right) + 
\frac\ga 2
\log\left(\frac{1-{\rm e}^{-\,\te\,R_t}}{1-{\rm
  e}^{-\,\te\,x}}\right)
%\log\left((1-{\rm e}^{-\,\te\,R_t})/(1-{\rm
%  e}^{-\,\te\,x})\right)
+
\frac {\te^2\,\ga} {8}\, 
\int_0^t({\rm e}^{\,\te\,R_s}-1)^{-1}\,ds\\
&&\\
&&
\hskip4cm
-
\frac {\te^2\,\ga} {4}(\frac{\ga}{2}-1)\,
\int_0^t({\rm e}^{\,\te\,R_s}-1)^{-2}\,ds
\end{eqnarray*}

In view of Proposition \ref{occu2} we consider now the functional
(assume that the continuous function $\omega$ is such that the 
functional is well-defined) 
$$
\widehat I_t:=
\frac{\te^2}{2}\int_0^t
\left(\exp(\te\,\omega_s)-1\right)^{-1}ds.
$$
By absolute
continuity, 
$$
%\begin{eqnarray*}
\E_x^{\widehat Z}\left(\exp(-r\widehat I_t)\right)
=
\E_x^R\left(\exp(\widehat D_t-r\widehat I_t)\right),
$$
and, further, for $H_y:=\inf\{t\,:\, \om_t=y\}$ with $y>x$ 
$$
%\begin{eqnarray*}
\E_x^{\widehat Z}\left(\exp(-r\widehat I_{H_y})\right)
=
\E_x^R\left(\exp(\widehat D_{H_y}-r\widehat I_{H_y})\right).
$$
Proposition \ref{occu2} gives now
\begin{equation}
\label{ee1}
\E_{x^*}^{X}\left(\exp(-r\,H_{y^*})\right)
=
\E_x^R\left(\exp(\widehat D_{H_y}-r\widehat I_{H_y})\right),
\end{equation}
where $x^*=1-{\rm e}^{-\,\te\,x},$
$y^*=1-{\rm e}^{-\,\te\,y},$ and
$X$ is a Jacobi diffusion with parameters $\al, \be,$ and
$\ga=\al+\be\geq 2.$ The identity (\ref{ee1}) is equivalent with
\begin{equation}
\label{ee2}
\E_x^R\left(\exp(-\widehat I_{H_y}(p,q))\right)
=
\frac{y}{x}\left(\frac{1-{\rm e}^{-\,\te\,x}}
{1-{\rm e}^{-\,\te\,y}}\right)^{\ga/2}
\E_{x^*}^{X}\left(\exp(-r\,H_{y^*})\right),
\end{equation}
where
\begin{equation}
\label{pq}
p=\frac{\te^2}{4}\,(2r-\ga)\quad {\rm and}\quad
q=\frac{\te^2}{8}\,\ga\,(\ga-2).
\end{equation}
Letting $r=\al\be$ it is seen from (\ref{pq}),
after straightforward computations, that
$\al$ and $\be$ can be expressed as in (\ref{aa1})
and (\ref{bb1}), respectively. To conclude the proof it remains
to compute the Laplace transform on the right hand side of (\ref{ee2}) and
let $y\to +\infty.$ Because $x<y$ we have also $x^*<y^*$, and. hence,
\begin{equation}
\label{ee12} 
\E_{x^*}^{X}\left(\exp(-\al\be\,H_{y^*})\right)
= \frac{\psi^X_{\al\be}(x^*)}{\psi^X_{\al\be}(y^*)},
\end{equation}
where $\psi^X_{\al\be}$ is the unique (up to a multiple), positive 
increasing solution of
(\ref{gauss}) with $\ga=\al+\be.$ We remark that the uniqueness
follows from the fact that both boundaries, in this case, are 
entrance-not-exit. Moreover, because $H_{y^*}\to +\infty$ as $y^*\to
1$ (also as $y^*\to 0$) it must hold that 
$\psi^X_{\al\be}(y^*)\to +\infty$ as $y^*\to
1.$ 
Similarly, there exists an unique positive decreasing solution 
$\varphi^X_{\al\be}$ such that $\varphi^X_{\al\be}(y^*)\to +\infty$ as 
$y^*\to 0.$ From the general theory of differential equations it
follows that all other solutions of (\ref{gauss}) can be expressed as
linear combinations of $\psi^X_{\al\be}$ and $\varphi^X_{\al\be}.$
Recall (see Section \ref{sec1}) that the function 
$
x\mapsto F(\al,\be,\ga;x)
%\quad{\rm and}\quad   
%x\mapsto G(\al,\be,\ga;x):=F(\al,\be,\al+\be+1-\ga;x)
$
is a solution of (\ref{gauss}). Now $\ga=\al+\be$ and from
Abramowitz and Stegun \cite{abramowitzstegun70} 15.3.10 p. 559
\begin{equation}
\label{ee22} 
F(\al,\be,\al +\be;x)\simeq -\frac{\Ga(\al+\be)}
{\Ga(\al)\Ga(\be)}\,\log(1-x),\qquad {\rm as}
\quad x\uparrow 1,
\end{equation}
and, because $F(\al,\be,\al +\be;0)=1,$ it follows that 
$x\mapsto F(\al,\be,\al +\be;x) $ is increasing, and, consequently,  
$$
\psi^X_{\al\be}= F(\al,\be,\al +\be;\cdot).
$$
Therefore, from (\ref{ee12}) and (\ref{ee2}), 
%\begin{equation}
%\label{ee2}
$$
\E_x^R\left(\exp(-\widehat I_{H_y}(p,q))\right)
=
\frac{y}{x}\left(\frac{1-{\rm e}^{-\,\te\,x}}
{1-{\rm e}^{-\,\te\,y}}\right)^{\ga/2}
\frac{F(\al,\be,\al +\be;x^*)}
{F(\al,\be,\al +\be;y^*)}.
%\end{equation}
$$ 
Letting here $y\to +\infty$  and using (\ref{ee22}) proves the claim.
\end{proof}

\begin{remark} {\rm
The decreasing solution $\varphi^X_{\al\be}$
is given by 
$$
\varphi^X_{\al\be}(x)=\widehat F(\al,\be,\al +\be;x)
:=F(\al,\be,1;1-x).
$$
From the definition of $F$ and the facts that, in this case,
$\al\be>0$ and 
$$
(\al+k)(\be+k)=\al\be+k(\al+\be)+k^2,
$$ 
it is seen
that $\widehat F$  has the desired properties, i.e.,  
$$
\lim_{x\to 0}\widehat F(\al,\be,\al +\be;x)   =+\infty
\quad
{\rm and}
\quad  
\widehat F(\al,\be,\al +\be;1)   =1,
$$
implying that $\widehat F$ is decreasing. 
}
\end{remark} 

\begin{example}
{\rm  In the particular case $p=q=4\rho$, and $\te=2$, we have
$$
\alpha=\beta=\frac 12 +\frac 12\,\sqrt{1+8\rho}
$$
and, therefore 
$$
\E_x\exp\Bigl(-\int_0^\infty
 \frac{\rho ds}{\sh^2(R_s)}\Bigr)
= 
\frac{(1-{\rm e}^{-2\, x})^{2\,\al}\, \Ga^2(\al)}{2\,x\,\Ga(2\al)}\
F\big(\al,\al,2\al;1-{\rm e}^{-2\, x}\big)
$$

}
\end{example}
\medskip
\begin{example}{\rm Choosing $q=0$ yields
$$
\alpha= 1+ \frac{\rm i}\te\, \sqrt{2p},\quad 
\beta=1- \frac {\rm i}\te\, \sqrt{2p}=:\bar \al .
$$
Hence, $\al +\beta=2$ and 
\begin{eqnarray*}
&&
\E_x\exp\Bigl(-\int_0^\infty\frac{p}{\exp(\te R_s)-1}ds\Bigr)\\
&&\hskip3cm
=
\frac{1-{\rm e}^{-\te x}}{x\,\te}
\, \Ga(\al)\,\Ga(\bar\al)\,
F\big(\al,\bar\al,2;1-{\rm e}^{-\te x}\big)\\
&&\hskip3cm
=
\frac{1-{\rm e}^{-\te x}}{x\,\te}
\, \frac{\pi \sqrt{2p}}{\te\,\sinh(\pi \sqrt{2p}/\te)}\,
F\big(\al,\bar\al,2;1-{\rm e}^{-\te x}\big),
\end{eqnarray*}
where the formula (see \cite{abramowitzstegun70} formulae 6.1.28 and 6.1.29.)
$$
\Gamma(1+{\rm i}\,y)\,\Gamma(1-{\rm i}\,y)=\frac{\pi y}{\sinh(\pi y)}
$$
is used.
Letting $x\to 0$ yields
\begin{equation}
\label{dy11}
\E_0\exp\Bigl(-\int_0^\infty\frac{p}{\exp(\te R_s)-1}ds\Bigr)
=
\frac{\pi \sqrt{2p}}{\te\,\sinh(\pi \sqrt{2p}/\te)}.
\end{equation}
The term on the right hand side is the Laplace transform of
$H_{\pi/\te}(R)$ where the BES(3) process $R$  is started at 0, i.e.,
(\ref{dy11}) is equivalent with the following identity due to 
Donati--Martin and Yor \cite{donatimartinyor97} p. 1044:
\begin{equation}
\label{dy12}
\int_0^\infty\frac{1}{\exp(\te R_s)-1}ds
\quad{\mathop=^{\rm{(d)}}}\quad
H_{\pi/\te}(R).
\end{equation}
The derivation of (\ref{dy12}) in \cite{donatimartinyor97} is very
different than the one presented above, and is based on 
a formula for the Laplace transform of an
integral functional of a two-dimensional Bessel process.
}
\end{example}
\bigskip
\noindent
{\bf Acknowledgement.} We thank Marc Yor for comments and, in
particular, for several references for Jacobi diffusions.

\bibliographystyle{plain}
%\bibliography{andrei}
\bibliography{yor1}
\end{document}